\documentclass[11pt]{article}
\usepackage{amsmath,amsthm}
\usepackage{amssymb}
\usepackage{latexsym,boxedminipage}
\usepackage{epsfig}

\renewcommand{\ell}{l}

\def\G#1#2{G_{#1}^{(#2)}}
\def\SS#1#2{\S_{#1}^{(#2)}}
\newcommand{\pib}{\bar \pi}

\newcommand{\la}{\lambda}

\renewcommand{\qed}{\hfill{$\Box$}}
\newcommand{\Z}{{\mathbb Z}}
\newcommand{\R}{{\mathbb R}}
\renewcommand{\P}{{\mathcal P}}
\renewcommand{\S}{\mathfrak S}

\def\1{\mathrm{Id}}

\newtheorem{lemma}{Lemma}

\newtheorem{theorem}{Theorem}

\title{Symmetrically Constrained Compositions}

\author{Matthias Beck
\thanks{Research supported in part by NSF grant DMS-0810105}\\
Department of Mathematics\\
San Francisco State University\\
San Francisco, CA 94132\\
{\tt beck@math.sfsu.edu}
\and
Ira M.~Gessel\\
Department of Mathematics \\
Brandeis University \\
Waltham, MA 02454--9110 \\
{\tt gessel@brandeis.edu}
\and
Sunyoung Lee
\thanks{Research supported in part by NSF grant DMS-0300034}\\
{Computer Science}\\
{N. C. State University }\\
{ Raleigh, NC 27695}\\
{\tt slee7@unity.ncsu.edu} 
\and Carla D. Savage
\thanks{Research supported in part by NSF grants DMS-0300034 and INT-0230800}\\
{Computer Science}\\
{N. C. State University }\\
{ Raleigh, NC 27695}\\
{\tt savage@csc.ncsu.edu} }

\begin{document}
\maketitle
\begin{center}
\emph{Dedicated to George Andrews on the occasion of his 70th birthday.}
\end{center}
\begin{abstract}
\noindent
Given integers $a_1, a_2, \dots, a_n$,
with $a_1 + a_2 + \dots + a_n \geq 1$,
a symmetrically constrained composition
$\la_1 + \la_2 + \dots + \la_n = M$ of $M$ into $n$ nonnegative parts is one
that satisfies each of the the $n!$ constraints
 $\left\{\,\sum_{i=1}^n a_i \lambda_{ \pi(i) } \geq 0 : \, 
\pi \in S_n \,\right\}$.
We show how to compute the generating function of these compositions,
combining methods from partition theory, permutation statistics, and
lattice-point enumeration.
\end{abstract}



\section{Introduction}

\subsection{Constrained Compositions}

This work was inspired by the ``constrained compositions'' introduced
by Andrews, Paule, and Riese in \cite{PAVII}.
We consider the problem of enumerating
{\em symmetrically constrained compositions}, that is, {compositions}
of an integer $M$ into $n$ nonnegative parts
\[
M = \la_1 + \la_2 + \cdots + \la_n = |\la| \, ,
\]
where the sequence $(\la_1, \la_2, \dots, \la_n)$ is
{constrained} to satisfy 
a {symmetric} system of linear inequalities.
 For example, the compositions
$\la_1+ \la_2+ \la_3$ of $M$
satisfying
\begin{equation}
\la_{\pi(1)} + \la_{\pi(2)} \geq \la_{\pi(3)}
\label{triangles}
\end{equation}
for every permutation $\pi \in S_3$,
are known as {\em integer-sided triangles}
of perimeter $M$
 \cite{andrews,PA2,JWW,stanley1}.
The number $\Delta_M$ of {\em incongruent} triangles of perimeter $M$ is given by
\begin{equation*}
\sum_{ M \ge 0 } \Delta_M \, q^M \   =  \  
\sum_{\substack{\la_1 \geq \la_2 \geq \la_3 \geq 0 \\[2pt]
\la_2 + \la_3 \geq \la_1 }} q^{|\la|}
 \  = \  
\frac{1}{(1-q^2)(1-q^3)(1-q^4)} \, .
\end{equation*}
(Note that, in contrast to  \cite{andrews,PA2,JWW,stanley1}, we allow
$\la_i=0$.)

However, $3+2+1$ and $2+3+1$ are different {\em compositions} (i.e., ordered partitions) of 6 and
counting the number  $\Delta^*_M$ of {\em ordered} solutions to \eqref{triangles} gives
\begin{equation}
\sum_{ M \ge 0 } \Delta^*_M \, q^M \  =  \ 
\sum_{\substack{\la_1 + \la_2 \geq \la_3  \\[2pt]
\la_1 + \la_3 \geq \la_2  \\[2pt]
\la_2 + \la_3 \geq \la_1}} q^{|\la|} \ 
 = \ 
\frac{1+2q^2 + 2q^4 + q^6}{(1-q^2)(1-q^3)(1-q^4)}  \ = \ 
\frac{1+q^3}{(1-q^2)^3} \, .
\label{ist}
\end{equation}
One could generalize this example in several ways.
For example, moving to $n$ dimensions, one could ask for the integer
sequences $(\la_1, \la_2, \dots, \la_n)$ satisfying
\begin{equation*}
\la_{\pi(1)}+\la_{\pi(2)}+\cdots+\la_{\pi(n-1)}\ge\la_{\pi(n)}
\end{equation*}
for {all} $n!$  permutations
$\pi \in S_n$.
Another generalization would be to study,
given positive integers $k,\ell,m$,
the integer
sequences $(\la_1, \la_2, \dots, \la_n)$ satisfying
\[
k\la_{\pi(1)} + \ell \la_{\pi(2)} \geq m\la_{\pi(3)}
\qquad \text{ for all } \pi \in S_n \,.
\]

Another example, considered in \cite{PAVII}, was inspired by a Putnam
exam problem \cite[Problem B3]{putnam-prob}: Enumerate all compositions of $M = \la_1 + \la_2$ into two parts satisfying
\[
2 \la_1 \geq \la_2 \qquad \text{ and } \qquad
2 \la_2 \geq \la_1 \, .
\]
It is shown that
\begin{equation}
\sum_{\substack{ 2 \la_1 \geq \la_2\\[2pt]
2 \la_2 \geq \la_1
}} x^{\la_1}y^{\la_2}
 \ = \ 
\frac{1+xy + x^2y^2}{(1-xy^2)(1-x^2y)} \, , 
\label{putnam}
\end{equation}
giving a complete parametrization of all solutions.

In \cite{PAVII}, Andrews, Paule, and Riese demonstrate the suitability of the
Omega package \cite{PA3} for experimenting with problems of this sort and
the power of MacMahon's partition analysis \cite{PA2} to prove some
elegant generalizations.

The goal of this paper is 
(1) to formulate a generalization of the symmetrically constrained compositions enumeration problem;
(2) to show how this problem is connected to permutation statistics;
(3) to show that the  permutation statistics approach gives, for many
cases, an effective computation method and, for certain cases, a way to
derive compact formulas; and
(4) to show that the insight provided by the geometry of lattice-point
enumeration aids in the handling of the most general case.

\subsection{The Symmetrically Constrained Compositions Enumeration Problem}

Let $\Z_{ \ge 0 }$ denote the set of nonnegative integers.
Fix integers $a_1, a_2, \dots, a_n$ (which may be negative).
 We are interested in enumerating 
compositions $\lambda = \left( \lambda_1, \lambda_2, \dots, \lambda_n \right) \in
 \Z_{ \ge 0 }^n$ that satisfy the $n!$ homogeneous linear constraints
\[
  a_1 \lambda_{ \pi(1) } + a_2 \lambda_{ \pi(2) } + \dots + a_n \lambda_{ \pi(n)
 } \ge 0
  \qquad \text{ for all } \pi \in S_n \, .
\]
Specifically, we are interested in computing the generating functions
\[
  F \left( z_1, z_2, \dots, z_n \right) := \sum_{ \lambda } z_1^{ \lambda_1 } z_
2^{ \lambda_2 } \cdots z_n^{ \lambda_n }
\]
and
\[
  F(q) := F(q,q,\dots,q) = \sum_{ \lambda } q^{ \lambda_1 + \lambda_2 + \dots +
\lambda_n }, 
\]
by exploiting the symmetry of the constraints.
Note that
because of the symmetry, there is no loss of generality in assuming
that
\[
a_1 \leq a_2 \leq \dots \leq a_n \, ,
\]
which we will do from now on.

In Section 2, we show how to solve the enumeration problem
when
$\sum_{i=1}^{n} a_i = 1$.
In certain special cases, we show that permutation statistics can be used
to derive elegant formulas.  
We note that even this simple case is
 difficult for general purpose software
 like the Omega Package \cite{PA3},
Xin's improvement of Omega \cite{xin}, and LattE macchiato \cite{LattE,Koeppe},
designed to enumerate solutions to linear Diophantine equations and
inequalities.  In Section 3 we solve the general problem.
We close this section with some notation and background on
permutation statistics.

\subsection{Permutation Statistics}

Throughout the paper, the following notation is used: 
 $[\, n \,] _q = (1-q^n)/(1-q)$;
 $[\, n \,] _q ! = \prod_{i=1}^n [\, i \,] _q$; and
  $(a;q)_n =
\prod_{i=0}^{n-1} (1-aq^i)$.

 For a permutation $\pi =
\pi(1)\pi(2) \cdots \pi(n)$ of $[n] := \{ 1, 2, \dots, n \}$,
the {\em descent set of $\pi$} is
\[
D_{\pi} = \{\,j : \, \pi(j) > \pi(j+1) \,\} \, .
\]
 The statistic $\mathrm{des}(\pi)=|D(\pi)|$ is
the number of descents of $\pi$ and the {\em major index} of $\pi$
is the sum of the descent positions: $\mathrm{maj}(\pi)=\sum_{i \in
D(\pi)}i$.
It is well known that
\begin{equation}
\sum_{\pi \in S_{n}}
q^{\mathrm{maj}(\pi)} \ = \
\prod_{i=1}^{n}\frac{1-q^i}{1-q} = [\,n\,]_q!
\label{eq:majdist}
\end{equation}
(see, e.g., \cite{stanley1}).
The joint distribution of $\mathrm{des}(\pi)$ and $\mathrm{maj}(\pi)$
over the set $S_n$ of all permutations of $[n]$ is given by Carlitz's
$q$-Eulerian polynomial
\cite{carlitz1, carlitz2}:
\[
C_n(x,q) = \sum_{\pi \in S_{n}}
x^{\mathrm{des}(\pi)}q^{\mathrm{maj}(\pi)} \ = \
\prod_{i=0}^{n}(1-xq^i) \sum_{j=1}^{\infty} [\, j\,] _q^n \, x^{j-1} .
\]
Applying the definition of  $[\, j\,] _q$ and the binomial theorem, we can
rewrite this as
\begin{equation}
C_n(x,q) = 
\frac{(x;q)_{n+1}}{(1-q)^n} \sum_{i=0}^n {n \choose i}
\frac{(-q)^i}{1-q^ix} \, . 
\label{eq:carlitz2}
\end{equation}
So, for example,
\begin{align}
C_1(x,q) & =  1 \nonumber\\
C_2(x,q) & =  1+ xq \label{C123} \\
C_3(x,q) & =  1+2xq+2xq^2+x^2q^3 . \nonumber
\end{align}

If we take the limit as $x\to q^{-n}$ in \eqref{eq:carlitz2} all terms except $i=n$ in the sum are canceled by $(q^{-n};q)_{n+1}=0$ in the numerator, so
\begin{equation}
\label{eq:limit}
C_n(q^{-n},q) = \frac{(-q)^n}{(1-q)^n}\lim_{x\to q^{-n}} \frac{(x;q)_{n+1}}{1-q^n x}
=\frac{(-q)^n}{(1-q)^n}\lim_{x\to q^{-n}}{(x;q)_{n}}
=\frac{(-q)^n}{(1-q)^n}{(q^{-n};q)_{n}} \, .
\end{equation}

Finally, for  $i \leq n-1$,
let $S_{n}^{(i)}$ be the set of permutations of $[n]$
that have no
descent in positions $\{n-i,n-i+1, ..., n-1 \}$.
Let 
\begin{align*}
C_{n}^{(i)}(x,q) & : =  \sum_{\pi \in S_{n}^{(i)}}
x^{\mathrm{des}(\pi)}q^{\mathrm{maj}(\pi)}. \label{Fnidef}
\end{align*}
In \cite{permstats}, it is shown that
\begin{equation*}
C_{n}^{(i)}(x,q)=\frac{C_{n}(x,q)}{(xq^{n-i};q)_i}-\sum_{k=1}^i
{n\choose k}xq^{n-k}
\frac{C_{n-k}(x,q)}{(xq^{n-i};q)_{i-k+1}}
\label{eq:sc}
\end{equation*}
so, in particular,
\begin{align}
C_{n}^{(1)}(x,q) & =  \frac{C_n(x,q)-nxq^{n-1}C_{n-1}(x,q)}
{1-xq^{n-1}} \, .
\label{C1}
\end{align}

\section{Symmetrically Constrained Compositions when
$\sum a_i = 1$}

\subsection{The Main Theorem}

\begin{theorem}\label{mainthm}
Given integers $a_1 \leq a_2 \leq \dots \leq a_n$ satisfying
$\sum_{i=1}^n a_i = 1$, the generating function for those $\la \in \Z^n_{ \ge 0 } $
satisfying
\[
\sum_{ j=1 }^n a_j \lambda_{ \pi(j) } \ge 0 \ \ \qquad \text{ for all } \pi \in S_n \,
\]
is
\begin{align*}
F(z_1,z_2, \dots, z_n) & = 
\sum_{\pi \in S_n}
\frac{\prod_{j\in D_{\pi}}\left(z_{\pi(1)}^{b_{1,j}}z_{\pi(2)}^{b_{2,j}} \cdots 
z_{\pi(n)}^{b_{n,j}}\right)}
{\prod_{j=1}^n  \left(1-z_{\pi(1)}^{b_{1,j}}z_{\pi(2)}^{b_{2,j}}
 \cdots z_{\pi(n)}^{b_{n,j}}\right) }
\end{align*}
where
\[
b_{i,j} = \left \{
\begin{array}{ll}
1 & \mbox{if $j = n$,}\\
-(a_1+ \dots + a_j) & \mbox{if $n \geq i > j \geq 1$,}\\
1 -(a_1+ \dots + a_j) & \mbox{if $1 \leq i \leq j < n$}.
\end{array}
\right.
\]
In particular, setting $z_1 = \cdots = z_n = q$ yields
\begin{align*}
F(q) & = 
\frac{\sum_{\pi \in S_n}\prod_{j \in D_{\pi}}q^{j-n\sum_{i=1}^j a_i}}
{(1-q^n)\prod_{j=1}^{n-1}\left(1-q^{j-n\sum_{i=1}^j a_i}\right)} \, .
\end{align*}
\end{theorem}
\noindent
{\bf Proof.}
To simplify notation, let
\[
F(z) = 
F(z_1,z_2, \dots, z_n) \, . 
\]
For $b \in \Z^n$, let
\[
z^b = z_1^{b_1} z_2^{b_2} \cdots  z_n^{b_n}
\]
and for $\pi \in S_n$, let
\[
z_{\pi} = (z_{\pi(1)}, z_{\pi(2)}, \dots, z_{\pi(n)}) \, .
\]
With
\[
  L := \biggl\{ \la  \in \Z_{ \ge 0 }^n : \, 
\sum_{ j=1 }^n a_j \la_{ \pi(j) } \ge 0  \text{ for all } \pi \in S_n \,
\biggr\}
\]
we have
\[
  F \left( z \right) = \sum_{ \lambda \in L } z^{\la}.
\]
Now we use the
standard method of partitioning the elements of $L$  into classes
$L_{\pi}$ indexed
by permutations $\pi \in S_n$:
\begin{eqnarray*}
L_{\pi} = \Bigl\{\, \la \in \Z^n & : & \la_{\pi(1)} \geq \la_{\pi(2)} \geq \dots \geq \la_{
\pi(n)},\\
 && \sum_{i=1}^{n} a_i \la_{\sigma(i)} \geq 0 \ \text{ for all $\sigma\in S_n$, and}\\
&& \la_{\pi(i)} > \la_{\pi(i+1)} \ {\rm if} \ i \in D_{\pi}\,\Bigr\} \, .
\end{eqnarray*}
Since
the last condition guarantees that no $\la$ is in more than
one class, $L$ is the disjoint union
\[
L = \bigcup_{\pi \in S_n} L_{\pi} \, .
\]
Our goal now simplifies to computing
\[
  F_\pi \left( z \right) := \sum_{ \lambda \in L_\pi} z^\lambda ,
\]
because $F \left( z \right) = \sum_{ \pi \in S_n } F_\pi \left( z \right)$.
In $L_{\pi}$, since 
\[
 \la_{ \pi(1) } \ge \la_{ \pi(2) } \ge \dots \ge \la_{ \pi(n) } \ge 0,
\]
and since, by our assumption, $a_1 \leq a_2 \leq \dots \leq a_n$,
the $n!$ constraints 
\[
 a_1 \lambda_{ \sigma(1) } + a_2 \lambda_{ \sigma(2) } + \dots + a_n \lambda_{ \sigma(n) } \ge 0
\qquad \text{ for all } \sigma \in S_n \,
\]
 are all implied by the single constraint
\[
  a_1 \la_{ \pi(1) } + a_2 \lambda_{ \pi(2) } + \dots + a_n \lambda_{ \pi(n) } \ge
 0 \, ,
\]
so that we get the more compact description
\[
  L_\pi = \left\{ \la \in \Z^n : \,
  \begin{array}{l}
    \la_{ \pi(1) } \ge \la_{ \pi(2) } \ge \dots \ge \la_{ \pi(n) } \ge 0 \ \text{ and
} \ \la_{ \pi(j) } > \la_{ \pi(j+1) } \text{ if } j \in D_\pi \\
    a_1 \la_{ \pi(1) } + a_2 \la_{ \pi(2) } + \dots + a_n \la_{ \pi(n) } \ge 0
  \end{array}
  \right\}.
\]
But this means that all $L_\pi$ look similar, except for the strict inequalities
 determined by $D_\pi$. More precisely, if we let
\[
  \widetilde L_\pi := \left\{ \la \in L_\1 : \, \la_j > \la_{ j+1 } \text{ if } j \in
D_\pi \right\}
\]
and $G_\pi (z) := \sum_{ \lambda \in \widetilde
L_\pi } z^\lambda$, then
\[
  F_\pi (z) = G_\pi \left( z_{\pi} \right) .
\]
So it remains to find $G_\pi( z) $, the generating function for
\[
  \widetilde L_\pi = \left\{ \la \in \Z^n : \,
  \begin{array}{l}
    \la_1 \ge \la_2 \ge \dots \ge \la_n \ge 0 \ \text{ and } \ \la_j > \la_{ j+1 } \text{
if } j \in D_\pi \\
    a_1 \la_1 + a_2 \la_2 + \dots + a_n \la_n \ge 0
  \end{array}
  \right\} ,
\]
for a given $\pi \in S_n$.

The constraints of $\widetilde L_\pi$, for $\la \in \Z^n_{\geq 0}$,  are given by the system
\begin{equation}\label{constraintmatrix1}
  \left[
    \begin{array}{cccccccccccccccccccc}
      1 & -1 & \\
        & 1 & -1 & \\
        &   &    & \ddots \\
        &   &    &        & 1 & -1 \\
     a_1&a_2& a_3& \cdots &a_{ n-1 }& a_n
    \end{array}
  \right]
  \ \la \ \ge \
  \left[
    \begin{array}{c}
      e_1 \\ e_2 \\ \vdots \\ e_{ n-1 } \\ e_n
    \end{array}
  \right] ,
\end{equation}
where
\[
  e_j =
  \begin{cases}
    0 & \text{ if } j \notin D_\pi \, , \\
    1 & \text{ if } j \in D_\pi \, .
  \end{cases}
\]
We make use of the following lemma, a well known result in
lattice-point enumeration.
This version
was formulated in \cite{Cmatrix2} for easy application to partition and 
composition
enumeration problems.
\begin{lemma}
Let $C = [c_{ i,j }]$ be an $n \times n$ matrix of integers such that $C^{-1} = B
= [b_{i,j}]$ exists and $b_{ i,j }$ are  all nonnegative integers. Let
$e_1, \dots, e_n$ be nonnegative integer constants.
For each $1 \leq i \leq n$, let
$c_i$ be the constraint
\[
c_{i,1}\lambda_1 +
c_{i,2}\lambda_2 + \dots +
c_{i,n}\lambda_n \geq e_i \, .
\]
Let $S_C$ be the set of nonnegative integer sequences $\la=
(\la_1,\la_2, \dots,\la_n)$ satisfying the constraints $c_i$ for
all $i$, $1 \leq i \leq n$. Then the generating function for $S_C$ is:
\[
F_C(x_1,x_2, \dots, x_n) =
\sum_{\lambda \in S_C}x_1^{\la_1} x_2^{\la_2} \cdots x_n^{\la_n} =
\frac{\prod_{j=1}^{n}\left(x_1^{b_{1,j}}x_2^{b_{2,j}} \cdots x_n^{b_{n,j}}\right)^{e_j}}
{\prod_{j=1}^n  \left(1-x_1^{b_{1,j}}x_2^{b_{2,j}}
 \cdots x_n^{b_{n,j}}\right) } \, .
\]
\label{Cmatrix}
\end{lemma}
Now let $C$ be the matrix on the left side of (\ref{constraintmatrix1}).
Then $\det(C) = a_1+\cdots + a_n=1$, so $C$ is invertible and
$B=C^{-1}$ has all integer entries:
\[
b_{i,j} = \left \{
\begin{array}{ll}
1 & \mbox{if $j = n$,}\\
-(a_1+ \dots + a_j) & \mbox{if $n \geq i > j \geq 1$,}\\
1 -(a_1+ \dots + a_j) & \mbox{if $1 \leq i \leq j < n$}.
\end{array}
\right.
\]
If, in addition,   $a_1 + \dots + a_j \leq
 0$  for
$1 \leq j \leq n-1$, the integer entries of $B=C^{-1}$ are all nonnegative
and Lemma 1 gives the generating function $G_{\pi}(z)$ and
the theorem follows.

To complete the proof, we show that
if $a_1 + \dots + a_n = 1$ and $a_1 \leq a_2 \leq \dots \leq a_n$,
then  for $1 \leq j \leq n-1$ we have
  $a_1 + \dots + a_j \leq 0$.

Let $j$ be the smallest index satisfying $1 \leq j \leq n-1$ and
$a_1 + \dots + a_j \leq 0$, but $a_1 + \dots + a_{j+1} > 0$.
Then $a_{j+1} > -(a_1 + \dots + a_j) \geq 0$.
Thus
\[
1 \leq a_{j+1} \leq \dots \leq  a_n \, ,
\]
so
\[
1 = a_1 + \dots + a_n \geq a_1 + \dots + a_{j+1} + n-j-1 > n-j-1 \, .
\]
So $j=n-1$ and therefore  $a_1 + \dots + a_j \leq 0$ for
$1 \leq j \leq n-1$.
\qed

In Section 2.3 we derive an algorithm based on Theorem 1
 for efficient computation of $F(q)$, given the $a_i$.
In the next section, we give examples of how to combine Theorem 1 with results on permutation statistics to derive formulas for $F(q)$
in special cases.

\subsection{Applications}

\noindent
{\bf Example 1}
Given positive integers $b$ and $n \geq 2$,
 let $L$ be the set of nonnegative integer sequences
$\la$
satisfying 
\[
(nb-b+1)\la_{\pi(n)} \geq b(\la_{\pi(1)} + \cdots + \la_{\pi(n-1)})
\qquad \text{ for all } \pi \in S_n \,.
\]
The case $n=2$, $b=1$ is the Putnam problem (\ref{putnam}).
Here $a=[-b,-b, \dots, -b, nb-b+1]$, so by Theorem 1,
\begin{equation*}
F(q) \   = \  
\frac{\sum_{\pi \in S_n}\prod_{j \in D_{\pi}}q^{j+jbn}}
{(1-q^n)\prod_{j=1}^{n-1}(1-q^{j+jbn})}\\
  \ =  \ 
\frac{\sum_{\pi \in S_n} (q^{1+bn})^{\mathrm{maj}(\pi)}}
{(1-q^n)\prod_{j=1}^{n-1}(1-q^{j+jbn})} \, .
\end{equation*}
By (\ref{eq:majdist}), the numerator is just $[\,n\,]_{q^{1+bn}}!$ and simplifying gives
\[
F(q) =  
\frac{1-q^{n(nb+1)}}
{(1-q^n)(1-q^{nb+1})^n} \, .
\]
This generating function was discovered by Andrews, Paule, and Riese and
a complete parametrization was proved in \cite{PAVII} using partition analysis.

\noindent
{\bf Example 2}
Given positive integers $b$ and $n \geq 2$, let $L$
be the set of nonnegative integer sequences
$\la$ satisfying 
\[
b(\la_{\pi(2)} + \cdots + \la_{\pi(n-1)}) \geq
(nb-b-1)\la_{\pi(1)}
\qquad \text{ for all } \pi \in S_n \,.
\]
The case $n=3$, $b=1$ is the integer-sided triangle problem
(\ref{ist}) and the case $n=2$, $b=2$ is the Putnam problem (\ref{putnam}).
Here $a = [-(nb-b-1),b,b, \dots, b]$,
 so by Theorem 1,
\begin{equation*}
F(q) \   =  \ 
\frac{\sum_{\pi \in S_n}\prod_{j \in D_{\pi}}q^{(bn-1)(n-j)}}
{(1-q^n)\prod_{j=1}^{n-1}(1-q^{(bn-1)(n-j)})}\\
  \ =  \ 
\frac{\sum_{\pi \in S_n} (q^{1-bn})^{\mathrm{maj}(\pi)}
(q^{n(bn-1)})^{\mathrm{des}(\pi)}}
{(1-q^n)\prod_{j=1}^{n-1}(1-q^{(bn-1)(n-j)})} \, .
\end{equation*}
By (\ref{eq:limit}), the numerator is
\begin{equation*}
C_n(q^{n(bn-1)}, q^{1-bn})
 =  \frac{(q^{n(bn-1)};q^{1-bn})_{n} (-q)^{(1-bn)n}}{(1-q^{1-bn})^n} \, ; 
\end{equation*}
Simplifying further and dividing by the denominator gives
\[
F(q) =
\frac{1-q^{n(nb-1)}}
{(1-q^n)(1-q^{nb-1})^n} \, .
\]
This generating function was also originally proved by Andrews, 
Paule, and Riese in \cite{PAVII}.

\noindent
{\bf Example 3}
Given positive integers $b$ and $n \geq 2$, let $L$ be the set 
of
nonnegative integer sequences
$\la = (\la_1, \dots,  \la_n)$ satisfying the constraints
\[
(b+1)\la_{\pi(n)} \geq b \la_{\pi(1)}
\qquad \text{ for all } \pi \in S_n \, .
\]
The case $n=2$, $b=1$ is the Putnam problem (\ref{putnam}).
Here $a=[-b,0,0, \dots, 0,b+1]$
 so by Theorem 1,
\begin{equation*}
F(q) \  = \ 
\frac{\sum_{\pi \in S_n}\prod_{j \in D_{\pi}}q^{j+bn}}
{(1-q^n)\prod_{j=1}^{n-1}(1-q^{j+bn})}\\
 \ = \ 
\frac{\sum_{\pi \in S_n} q^{\mathrm{maj}(\pi)}
(q^{bn})^{\mathrm{des}(\pi)}}
{(1-q^n)\prod_{j=1}^{n-1}(1-q^{j+bn})} \, .
\end{equation*}
By (\ref{eq:carlitz2}), the numerator is
\begin{align*}
C_n(q^{bn}, q) & = 
\frac{(q^{bn};q)_{n+1}}{(1-q)^n} \sum_{i=0}^n {n \choose i}
\frac{(-q^{i})}{1-q^{bn+i}}  \, .
\end{align*}
Combining with the denominator and simplifying gives
\[
F(q) = \frac{(1-q^{bn})(1-q^{bn+n})}
{(1-q^n)(1-q)^n} \sum_{i=0}^n {n \choose i}
\frac{(-q^{i})}{1-q^{bn+i}} \, .
\]

\noindent
{\bf Example 4}
Given  positive integers $k \leq \ell$, and $n \geq 3$,
let $m=k+\ell-1$ and
 let $L$ be the set of
nonnegative integer sequences
$\la$
satisfying
\[
k \la_{\pi(n-1)} + \ell \la_{\pi(n)} \geq m \la_{\pi(1)}
\qquad \text{ for all } \pi \in S_n \,.
\]
The case $n=3$ and $k=\ell=1$ is the integer-sided triangles (\ref{ist}).
Here $a=[-m, 0,0,\dots, 0, k,\ell]$, so by Theorem 1,
\begin{align*}
F(q) & = 
\frac{\sum_{\pi \in S_n}\prod_{j \in D_{\pi}, j \not = n-1}q^{j+mn}
\prod_{j \in D_{\pi}, j = n-1}q^{j+mn -nk}}
{(1-q^n)(1-q^{n\ell -1})\prod_{j=1}^{n-2}(1-q^{j+bn})} \, .
\end{align*}
Recall from (\ref{C1}) that $C_n^{(1)}$ is the joint distribution
of des and maj over all permutations with no descent in position $n-1$.
Then in $F(q)$, we can
split the sum over $\pi \in S_n$ into two sums according to
whether or not $i \in D_{\pi}$.
We get that the numerator can be written as:
\[
C_n^{(1)}(q^{nm},q) + q^{-nk} (C_n(q^{nm},q) - C_n^{(1)}(q^{nm},q)) \, .
\]
Using (\ref{C1}) and combining with the denominator gives
(eventually)
\begin{align*}
F(q) & = 
\frac{C_n(q^{nm},q)(1-q^{n\ell-1})
- C_{n-1}(q^{nm},q) nq^{nm+n-1}(1-q^{-nk})}
{(1-q^n)(1-q^{n\ell-1})(q^{nm+1};q)_{n-1}} \, .
\label{fnq}
\end{align*}
Details appear in \cite{permstats}.

\subsection{Efficient Enumeration of Symmetrically Constrained Compositions}
We can compute the generating function $F(q)$ for compositions satisfying the $n!$ constraints
\[
\sum_{i=1}^n a_i \la_{\pi(i)} \geq 0
\qquad \text{ for all } \pi \in S_n \,
\]
via Theorem 1.
The denominator is given explicitly, but the numerator is a sum of $n!$
terms.  However, regardless of the values of the $a_i$, the numerator
of $F(q)$, when simplified, is a polynomial with at most $2^{n-1}$
terms (one for each possible descent set).

Let $u_1, u_2,\dots$ be arbitrary and define polynomials $G_n$ by
\begin{equation*}
G_n=\sum_{\pi\in S_n} \prod_{i\in D_{\pi} } u_i \, .
\end{equation*}
We can compute $G_n$ in the following way:
Let
\begin{equation*}
\G{n}{i} = \sum_{\pi\in\SS ni} \prod_{i\in D_{\pi} } u_i \, ,
\end{equation*}
where $\SS ni$ is the set of all  permutations $\pi\in S_n$ that end with $i$.
Then $G_n = \G{n+1}{n+1}$.

A permutation $\pi$ in $\SS ni$ can be obtained uniquely from some permutation
$\pib$ in $S_{n-1}$ by replacing each $j\ge i$ with $j+1$ and then appending $i
$ at the end. The descent set of $\pi$ will be the same as the descent set of 
$\pib$ if the last entry of $\pib$ is less than $i$ and the descent set of $\pi$ will be $D_{\pib}\cup\{n-1\}$ if 
the last entry of $\pib$ is greater than or equal to $i$. Thus we have the recurrence
\begin{equation*}
\G ni = \sum_{j=1}^{i-1} \G{n-1}j + u_{n-1}\sum_{j=i}^{n-1} \G{n-1}j
\end{equation*}
with the initial condition $\G 11=1$.
We can simplify this a bit to get {\bf ``Algorithm $G$''}:
\begin{align*}
\G ni =  \G n{i-1} +(1-u_{n-1})\G{n-1}{i-1} \quad\text{for $i>1$} \, ,
\end{align*}
with $\G n1 = u_{n-1}\sum_{j=1}^{n-1} \G{n-1}j$.

Now, to compute the numerator of $F(q)$ in Theorem 1 using Algorithm $G$,
simply set $u_i = q^{i-n(a_1 + \cdots + a_i)}$ for $1 \leq i < n$ and
compute  $\G{n+1}{n+1}$.

If we use dynamic programming to implement the recurrence of Algorithm $G$,
(e.g. ``option remember'' in Maple), then to compute $G_n = \G{n+1}{n+1}$,
at most
$O(n^2)$ polynomials are computed. However, we must consider the time required
to compute them.
In order to compute one of the $G_k^{(i)}$, essentially we
only need to add two polynomials.
It is fair to assume that the time is proportional
to the number of terms in the polynomials times the logarithm of the 
coefficient magnitude.  So, overall, the time (and number of terms) grows
roughly like $2^n$ in the dimension $n$, but logarithmically in the
coefficient size, which is considered polynomial time in fixed dimension.
In practice, we found that we could compute $F(q)$ for arbitrary $a$ with
$\sum a_i = 1$ within seconds for $n \leq 11 $, in about 10 seconds for
$n=12$  and
in less than a minute up to $n=15$, 
using a naive implementation in Maple on a tablet PC running Windows XP.

For comparison, there are existing software packages that,
when given a collection of linear inequalities, produce the
generating function for the integer points in the solution set.
These packages include the  Omega Package \cite{PA3},
Xin's speed-up of Omega \cite{xin}, and LattE macchiato \cite{LattE,Koeppe}.
We used these programs to compute symmetrically constrained compositions in $n$ 
dimensions, by giving as input the
$n!$ inequalities.  The computation became infeasible when $n \geq 4$
for the Omega package and Xin's program.  LattE was able to handle
examples for $n=5$ in under 10 seconds and $n=6$ in under an hour.

Thus exploiting the symmetry via Theorem 1 and Algorithm $G$ makes a huge difference in what we can compute.

\section{The General Case}

\subsection{A General Version of the Main Theorem}

We remove the requirement that $\sum a_i = 1$ and
enumerate compositions $\lambda = \left( \lambda_1, \lambda_2, \dots, \lambda_n \right) \in \Z_{ \ge
 0 }^n$ that satisfy the $n!$ constraints
\begin{equation}\label{symconstraintseq}
  a_1 \lambda_{ \pi(1) } + a_2 \lambda_{ \pi(2) } + \dots + a_n \lambda_{ \pi(n)
 } \ge 0
  \qquad \text{ for all } \pi \in S_n \, 
\end{equation}
via the generating function $F(z) = \sum_{\la} z^{\la}$.

\begin{theorem}\label{mainthmgeneralized}
Given integers $a_1 \leq a_2 \leq \dots \leq a_n$,
with $a_1 + a_2 + \dots + a_j \le 0$ for $1 \le j \le n-1$ and $a_1 + a_2 + \dots + a_n \geq 1$,
define the vectors $A_1, A_2, \dots, A_n 
\in \Z^n$ as the columns of the matrix
\[
  \left[
    \begin{array}{cccccccccccccccccccc}
      a_2 + \dots + a_n & a_3 + \dots + a_n & a_4 + \dots + a_n & \cdots & a_n &
 1 \\
      -a_1 & a_3 + \dots + a_n & a_4 + \dots + a_n & & a_n & 1 \\
      -a_1 & -a_1-a_2 & a_4 + \dots + a_n & & a_n & 1 \\
      -a_1 & -a_1-a_2 & -a_1-a_2-a_3 & & a_n & 1 \\
      \vdots & \vdots & \vdots & & \vdots & \vdots \\
      -a_1 & -a_1-a_2 & -a_1-a_2-a_3 & & a_n & 1 \\
      -a_1 & -a_1-a_2 & -a_1-a_2-a_3 & \cdots & -a_1-\dots-a_{ n-1 }  & 1
    \end{array}
  \right]
\]
and let
\[
  \P := \sum_{ j=1}^n  [0,1) A_j = \biggl\{\,\sum_{i=1}^n c_iA_j : \  0 \leq c_i < 1 \,\biggr\}
 \, .
\]
Then
\[
  F \left( z \right) = \sum_{p \in \P \cap \Z^n} \,   \sum_{ \pi \in S_n }
\frac{z_{\pi}^p \prod_{i \in D_{\pi},\, p_i = p_{i+1}}  z_{\pi}^{A_i}}
  {\prod_{  j=1 }^n  \left( 1 - z_{ \pi }^{A_i} \right)} \, ,
\]
where we take the product over all descent positions $i$ of $\pi$ for which the $i$th and the $(i+1)$st coordinate of $p$ are the same.
\end{theorem}

\noindent
If, for some $i$, $d$ divides every coordinate of $A_i$, we can
replace $A_i$ by $A_i/d$ in Theorem \ref{mainthmgeneralized} and thereby reduce the number of
lattice points in $\P$ by a factor of $d$.

\noindent
{\bf Proof.}
The start of our proof is similar to that of Theorem \ref{mainthm}, except that we find it advantageous to view the compositions satisfying 
\eqref{symconstraintseq} as integer points in the simplicial cone
\[
  K := \biggl\{\, x = \left( x_1, x_2, \dots, x_n \right) \in \R_{ \ge 0 }^n : \, 
\sum_{ j=1 }^n a_j x_{ \pi(j) } \ge 0\
 \text{ for all } \pi \in S_n \,
\biggr\} .
\]
A {\em simplicial cone} is a subset of $\R^n$ of the form
$\{ y \in \R^n \ | \ My \leq b \}$ where $M$ is a nonsingular real
matrix and $b \in \R^n$.
From this point of view,
\[
  F \left( z \right) = \sum_{ \lambda \in K \cap \Z^d } z^\lambda.
\]
The setup now continues in analogy with the proof of Theorem \ref{mainthm}. Like there, it suffices to study
\[
  \widetilde K_\pi := \left\{ x \in \R^n : \,
  \begin{array}{l}
    x_1 \ge x_2 \ge \dots \ge x_n \ge 0 \ \text{ and
} \ x_j > x_{ j+1 } \text{ if } j \in D_\pi \\
    a_1 x_1 + a_2 x_2 + \dots + a_n x_n \ge 0
  \end{array}
  \right\}
\]
and the associated generating function
$G_\pi (z) := \sum_{ \lambda \in \widetilde K_\pi \cap \Z^d } z^\lambda$; then
\[
  F \left( z \right) = \sum_{ \pi \in S_n } G_\pi \left( z_{ \pi } \right) .
\]
First, we study the cone $K_{\1}$.
The constraints of $ K_{\1}$ are given by the system
\begin{equation}\label{constraintmatrix}
  \left[
    \begin{array}{cccccccccccccccccccc}
      1 & -1 & \\
        & 1 & -1 & \\
        &   &    & \ddots \\
        &   &    &        & 1 & -1 \\
     a_1&a_2& a_3& \cdots &a_{ n-1 }& a_n
    \end{array}
  \right]
  \ x \ \ge \ 0 \, ,
\end{equation}
and the inverse of the matrix on the left of \eqref{constraintmatrix} is
\[
  \frac{ 1 }{ \sum_{ j=1 }^n a_j }
  \left[
    \begin{array}{cccccccccccccccccccc}
      a_2 + \dots + a_n & a_3 + \dots + a_n & a_4 + \dots + a_n & \cdots & a_n &
 1 \\
      -a_1 & a_3 + \dots + a_n & a_4 + \dots + a_n & & a_n & 1 \\
      -a_1 & -a_1-a_2 & a_4 + \dots + a_n & & a_n & 1 \\
      -a_1 & -a_1-a_2 & -a_1-a_2-a_3 & & a_n & 1 \\
      \vdots & \vdots & \vdots & & \vdots & \vdots \\
      -a_1 & -a_1-a_2 & -a_1-a_2-a_3 & & a_n & 1 \\
      -a_1 & -a_1-a_2 & -a_1-a_2-a_3 & \cdots & -a_1-\dots-a_{ n-1 }  & 1
    \end{array}
  \right] .
\]
The conditions on $a_1, a_2, \dots, a_n$ guarantee
that the inverse exists and that $K_{\1}$ is a cone in $\R^n_{\geq 0}$.
Thus the columns $A_1, A_2, \dots, A_n$ of this matrix form a set of 
generators of $K_{\1}$ and by an easy
tiling argument (see, e.g., \cite[Chapter 3]{beckbook})
\begin{equation}
K_{\1} \cap \Z_n = \biggl\{\, p + \sum_{j=1}^n c_jA_j : \, p \in \P, \ 
c \in \Z^n_{\geq 0} \,\biggr\} ;
\label{KId}
\end{equation}
in other words,
\[
  G_{\1} (z) =
\frac{
 \sum_{p \in \P \cap \Z^n} z^p}
 { \prod_{  j=1 }^n  \Bigl( 1 - z_{ \pi }^{A_j} \Bigr)} \, .
\]
Before turning to  $\widetilde K_\pi$, note that the generators $A_j$
have a special form: we have (writing $A_{i,j}$ for the $i$th entry of $A_j$)
\[
A_{i,j} =  \left \{
\begin{array}{ll}
1 & {\mbox{if $j=n$, else}}\\
-(a_1+a_2+ \ldots + a_j) & {\mbox{if $i>j$, else}}\\
a_{j+1}+a_{j+2}+ \ldots + a_n & 
\end{array}
\right . .
\]
Therefore, since $\sum_{ j=1 }^n a_j \geq 1 $, 
\[
A_{j,j} > A_{j+1,j} \text{ for } 1 \leq j < n
\qquad \text{ and } \qquad
A_{i,j} = A_{i+1,j} \text{ for } j \not = i.
\]
Thus, if $p \in K_{\1} \cap \Z^n$ satisfies $p_j = p_{j+1}$, then for
any $c \in \Z^n_{\geq 0}$, if
\[
r = (p+ A_j) + \sum_{i=1}^n c_i A_i
\]
then
\[
	r_j > r_{j+1}.
\]
Now, what about $\widetilde K_\pi$?
It contains all  points $y \in  K_{\1}$ except those $y$ with $y_i = y_{i+1}$
for some $i \in D_{\pi}$.
By (\ref{KId}), every $y \in  K_{\1}$ has a unique representation as
$y = p + \sum_{j=1}^n c_jA_j$ for some $c \in \Z^n_{\geq 0}$.
Thus by the previous paragraph, $y_i = y_{i+1}$ iff
both $p_i = p_{i+1}$ and $c_i=0$. 
Now, in the same way as in (\ref{KId}),
\begin{align*}
\widetilde K_\pi \cap \Z^n & =  \biggl\{\, p + \sum_{j=1}^n c_jA_j : \, p \in \P, \ c \in \Z^n_{\geq 0}, {\rm \  and  \ if \ } j \in D_{\pi} {\rm \  and \ }
p_j = p_{j+1} {\rm \ then \ } c_j > 0 \,\biggr\} \\
& = 
 \biggl\{\, p + \sum_{j=1}^n c_jA_j +\!\! \sum _{j \in D_{\pi}, \  p_j = p_{j+1}} \!\!\!\!\!\!\!\! A_j \, : \, p \in \P, \ c \in \Z^n_{\geq 0} \,\biggr\} 
\end{align*}
Thus
\[
G_{\pi}(z) =
 \sum_{p \in \P \cap \Z^n}
\frac{z^p \prod_{j \in D_{\pi}, \ p_j = p_{j+1} } z_{\pi}^{A_i}}
  {\prod_{  j=1 }^n  \left( 1 - z_{ \pi }^{A_i} \right)} \, . 
\]
\qed

\noindent
In the special case of Theorem 1, $\sum_{ j=1 }^n a_j = 1$ and the origin
is the only lattice point in $\P$.

\subsection{Efficient Computation for the General Case}

Give $a_1 \leq \dots \leq a_n$ with $\sum_{ j=1 }^n a_j \geq 1$, once we find the
generators $A_1, A_2, \dots, A_n$, and the lattice points in $\P$, we can
again use Algorithm $G$ to efficiently compute $F(q)$:
By Theorem~2, setting $z=(q,q, \dots, q)$,
\[
F(q) =
\sum_{p \in \P \cap \Z^n}  q^{|p|}
\frac{ \sum_{ \pi \in S_n }  \prod_{ i \in D_{\pi},\, p_i = p_{i+1} 
} q^{|A_i|}}
  {\prod_{  j=1 }^n  \left( 1 - q^{|A_i|} \right)},
\]
where $|x| = x_1 + \dots + x_n$ for an $n$-dimensional vector $x$.

The denominator is easy.  To find the numerator, for each
point $p \in \P \cap \Z^n$, set
\[
u_i = 
\left \{
\begin{array}{ll}
q^{|A_i|} & {\rm if \ } p_i = p_{i+1}\\
1 & {\rm otherwise}
\end{array}
\right.
\]
and then compute $\G{n+1}{n+1}$.

Now the running time also depends on $|\P \cap \Z^n|$.
This can grow linearly with the magnitude of the entries
(rather than the logarithm of the magnitude), even in fixed dimension.
However, when $|\P \cap \Z^n|$ is of moderate size, this computation method
can be quite effective.

{
\bibliography{symcc}
\bibliographystyle{plain}
}

\end{document}